\documentclass[12pt,twoside]{amsart}
\usepackage{amssymb}
\usepackage{amscd}

\title[Termination of $4$-fold canonical flips]
{Termination of 4-fold canonical flips}  
\author{Osamu Fujino} 
\date{2003/1/15}
\address{Research Institute for Mathematical Sciences\\ 
 Kyoto University, Kyoto 606-8502 Japan}
\email{fujino@kurims.kyoto-u.ac.jp}

\newcommand{\xdiscrep}[0]{{\operatorname{discrep}}}

\newcommand{\mult}[0]{{\operatorname{mult}}}
\newcommand{\rank}[0]{{\operatorname{rk}}}

\newtheorem{thm}{Theorem}
\newtheorem{lem}[thm]{Lemma}
\newtheorem{cor}[thm]{Corollary}

\theoremstyle{definition}
\newtheorem{defn}[thm]{Definition}
\newtheorem{rem}[thm]{Remark}

\theoremstyle{definition}

\newtheorem{say}[thm]{}

\theoremstyle{remark}
\newtheorem*{ack}{{\bf{Acknowledgments}}}         
 
\begin{document}
\bibliographystyle{amsalpha+}

\begin{abstract}
We prove the termination of 4-fold canonical flips. 
\end{abstract}

\maketitle

This paper is a supplement of \cite[Theorem 5-1-15]{kmmm} and 
\cite[Main Theorem 2.1]{matsuki}. 
We prove that a sequence of $4$-fold log flips for 
canonical pairs terminates after finitely many steps. 

Let us recall the definition of {\em{canonical flips}}, which is 
slightly different from the usual one (cf.~\cite[(2.11) Adjoint Diagram]
{S}). 

\begin{defn}[Canonical flip]\label{cflip}
Let $X$ be a normal projective variety and $B$ an 
effective $\mathbb Q$-divisor such that the pair 
$(X,B)$ is canonical, that is, $\xdiscrep(X,B)\geq 0$. 
Let $\phi:(X,B)\longrightarrow Z$ be a small contraction corresponding to a 
$(K_X+B)$-negative extremal face. 
If there exists a normal projective variety $X^+$ 
and a projective morphism
$\phi^+:X^+\longrightarrow Z$ such that
\begin{enumerate}
\item $\phi^+$ is small; 
\item $K_{X^+}+B^+$ is $\phi^+$-ample, where $B^+$ is 
the strict transform of $B$, 
\end{enumerate}
then we call $\phi^+$ the {\em{canonical flip}} or {\em{log flip}} 
of $\phi$.
We call the following diagram a {\em{flipping diagram}}: 
$$
\begin{matrix}
(X,B) & \dashrightarrow & {(X^+,B^+)} \\
{\phi\searrow} & \ &  {\swarrow \phi^+} \\
 \ & Z & \ & .
\end{matrix}
$$
\end{defn}

The following is the main theorem of this paper. 

\begin{thm}[Termination of $4$-fold canonical flips]\label{can}
Let $X$ be a normal projective $4$-fold and 
$B$ an effective $\mathbb Q$-divisor such 
that $(X,B)$ is canonical. 
Consider a sequence of log flips starting from $(X,B)=(X_0,B_0)${\em{:}}
$$
(X_0,B_0)\dasharrow (X_1,B_1)\dasharrow (X_2,B_2)\dasharrow \cdots, 
$$
where $\phi_i:X_i\longrightarrow Z_i$ is a contraction and 
${\phi_i}^{+}:{X_i}^{+}=X_{i+1}\longrightarrow Z_i$ is the log flip.  
Then this sequence terminates after finitely many steps. 
\end{thm} 

\begin{cor}[{\cite[Main Theorem 2.1]{matsuki}}]\label{flop}
Let $X$ be a projective $4$-fold with only terminal singularities 
and $D$ an effective $\mathbb Q$-Cartier $\mathbb Q$-divisor. 
Then any sequence of $D$-flops is finite. 
\end{cor}

We introduce a variant of {\em{difficulty}}. 
This is slightly different from \cite[Chapter 4]{FA}. 
It was inspired by \cite[4.14 Remark]{FA}. 
Note that the notion of difficulty was first introduced by Shokurov in 
\cite[(2.15) Definition]{S}.

\begin{defn}[A weighted version of difficulty]\label{diff}
Let $(X,B)$ be a pair with only canonical singularities, 
where $B=\sum _{j=1}^{l}b_jB^j$ 
with $0<b_1<\cdots<b_l\leq 1$ and $B^j$ is a reduced divisor for every $j$. 
We note that $B^j$ is not necessarily 
irreducible. 
We put $b_0=0$, and $S:=\sum _{j\geq 0}b_j\mathbb Z_{\geq 0}
\subset \mathbb Q$. 
Note that $S=0$ if $B=0$. 

We call a divisor $E$ over $X$ {\em{essential}} 
if  $E$ is exceptional over $X$ and is not obtained from blowing 
up the generic point of a subvariety 
$W\subset B\subset X$ such that $B$ and $X$ 
are generically smooth along $W$ (and thus only one 
of the irreducible components of 
$\sum_{j\geq 1}B^j$ contains $W$) and $\dim W=\dim X-2$. We set 
$$
d_{S,b}(X,B):=\sum _{\xi\in S, \xi\geq b}\sharp 
\{ E | E \ \text{is essential and} \ a(E,X,B)<1-\xi \}. 
$$ 
Then $d_{S,b_j}(X,B)$ is finite by Lemma \ref{finite} below. 
Note that $\xdiscrep (X,B)\geq 0$ since 
$(X,B)$ is canonical. 
\end{defn}

\begin{lem}[{\cite[(4.12.2.1)]{FA}}]\label{finite}
Let $(X,B)$ be a klt pair. Then 
$$
\sharp\{E| E\  {\text{is essential and 
}} a(E,X,B)<\min\{1,1+\xdiscrep (X,B)\}\}
$$ 
is finite. 
\end{lem}

\begin{defn}\label{type}
Let $\varphi:(X,B)\dasharrow (X^+,B^+)$ be a canonical flip. 
We say that this flip is {\em{of type $(\dim A, \dim A^+)$}}, 
where $A$ (resp.~$A^+$) is the exceptional locus of $\phi:X\longrightarrow Z$ 
(resp.~$\phi^+:X^+\longrightarrow Z$). We call $A$ (resp.~$A^+$) 
the {\em{flipping}} (resp.~{\em{flipped}}) 
locus of 
$\varphi$. 
When $\dim X=4$, the log flip is either of type $(1,2)$, $(2,2)$ or $(2,1)$ 
by \cite[Lemma 5-1-17]{kmmm}. 
\end{defn}

\begin{rem} 
In \cite[\S 5-1]{kmmm}, the variety is 
$\mathbb Q$-factorial and every flipping contraction is 
corresponding to a negative extremal ray. 
However, the above properties were 
not used in the proof of \cite[Lemma 5-1-17]{kmmm}.  
Note that \cite[Lemma 5-1-17]{kmmm} holds in more 
general setting. For the details, 
see the original article \cite{kmmm}. 
\end{rem}

\begin{lem}[{cf.~\cite[Lemma 6.21]{km}}]\label{keisan}
Let $\varphi:(X,\Delta)\dashrightarrow (X^+, \Delta^+)$ 
be a canonical flips of $n$-folds. 
Let $\Delta:=\sum \delta_i \Delta_i$ be the irreducible 
decomposition and $\Delta^{+}$ {\em{(}}resp.~$\Delta_i^{+}${\em{)}} 
the strict transform of $\Delta$ {\em{(}}resp.~$\Delta_i${\em{)}}. 
Let $F$ be an $(n-2)$-dimensional irreducible component 
of $A^+$, 
and $E_F$ the exceptional divisor obtained by blowing up $F$ 
near the generic point of $F$. 
Then $X^+$ is generically smooth along $F$ and 
$$
0\leq a(E_F, X,\Delta)<a(E_F, X^+, \Delta^+)=1-\sum \delta_i
\mult_F(\Delta_i^+), 
$$ 
where $\mult _F(\Delta_i^+)$ is the multiplicity 
of $\Delta_i^+$ along $F$. 
\end{lem}
\begin{proof}
The pair $(X^+,\Delta^+)$ is terminal near the generic point 
of $F$ by the negativity lemma. 
Therefore, $X^+$ is generically smooth 
along $F$, and the rest is an obvious computation. 
\end{proof}

\begin{say}
Let $X$ be a projective variety. 
Let $X^{an}$ be the underlying analytic space of 
$X$. 
Let $H_k^{BM}(X^{an})$ be the {\em{Borel-Moore}} 
homology. 
For the details, see \cite[19.1 Cycle Map]{fulton}. 
Then there exists a cycle map; 
$$
cl:A_k(X)\longrightarrow H_{2k}^{BM}(X^{an}), 
$$ 
where $A_k(X)$ is the group generated by 
rational equivalence classes of $k$-dimensional cycles on $X$. 
For the details about $A_k(X)$, see \cite[Chapter 1]{fulton}. 
We note that the cycle map 
$cl$ commutes with push-forward for proper morphisms, 
and with restriction to open subschemes. 
From now on, we omit the superscript $an$ for 
simplicity. 
\end{say}

The following is in the proof of \cite[Theorem 5-1-15]{kmmm}. 

\begin{lem}[Log flips of type $(2,1)$]\label{hom} 
When the $4$-dimensional log flip 
$$
\begin{matrix}
(X,B) & \dashrightarrow & {(X^+,B^+)} \\
{\phi\searrow} & \ &  {\swarrow \phi^+} \\
 \ & Z & \ & 
\end{matrix} 
$$ is of type $(2,1)$, 
we study the rank of the $\mathbb Z$-module $cl(A_2(X))$, 
where $cl(A_2(X))$ is the image of the cycle map 
$cl: A_2(X)\longrightarrow H_4^{BM}(X)$. 
By the following commutative diagrams{\em{;}} 
$$
\begin{CD}
A_2(A)@>>>A_2(X)@>>>A_2(X\setminus A)@>>>0\\
@VV{cl}V @VV{cl}V @VV{cl}V\\ 
H_4^{BM}(A)@>>>H_4^{BM}(X)@>>>H_4^{BM}(X\setminus A), 
\end{CD} 
$$
and 
$$
\begin{CD}
A_2(Z)@>{\sim}>>A_2(Z\setminus \phi(A))\\
@VV{cl}V @VV{cl}V\\ 
H_4^{BM}(Z)@>{\sim}>>H_4^{BM}(Z\setminus \phi(A)), 
\end{CD} 
$$
we have the surjective homomorphism{\em{;}} 
$$
cl(A_2(X))
\longrightarrow cl(A_2(Z)). 
$$ 
We note that $X\setminus A\simeq Z\setminus 
\phi(A)$ and \cite[p.371 (6) and Lemma 19.1.1]{fulton}. 
For any closed algebraic subvariety $V$ on $X$ of 
complex dimension $2$, 
$V$ is not numerically trivial since $X$ is projective. 
Therefore, $cl(V)\ne 0$ in $cl(A_2(X))$ 
{\em{(}}see \cite[Definition 19.1]{fulton}{\em{)}}. 
Thus the kernel of the surjection above is not zero. 
By the similar arguments, 
we obtain that 
$$
cl(A_2(X^+))
\simeq cl(A_2(Z)). 
$$
We note that $A^+$ is one-dimensional and 
$X^{+}\setminus A^{+}\simeq Z\setminus \phi^{+}(A^{+})$. 
Therefore, since the rank $\rank_{\mathbb Z} 
cl(A_2(X))$ is finite, we finally have the result
$$
\rank_{\mathbb Z} cl(A_2(X))>
\rank_{\mathbb Z} cl(A_2(X^+)). 
$$ 
\end{lem}

Let us start the proof of Theorem \ref{can}. 

\begin{proof}[Proof of {\em{Theorem \ref{can}}}]
Since $d_{S,b_l}(X,B)$ does not increase by 
log flips by the negativity lemma. 
Note that if $E$ is exceptional and $a(E,X_i, B_i)<1-b_l$, 
then $E$ is essential. If $B_{i}^l$, which is the strict transform of 
$B^l$ on $X_i$, 
contains $2$-dimensional flipped locus, 
then $d_{S,b_l}(X,B)$ decreases by easy computations 
(see Lemma \ref{keisan}). 
So, after finitely many flips, 
$B_{i}^{l}$ does not contain 
$2$-dimensional flipped locus. 
Thus, we can assume that $B_{i}^{l}$ does not 
contain $2$-dimensional flipped locus for 
every $i>0$. 
Let $\overline{B_{i}^{l}}$ be the normalization 
of $B_{i}^{l}$. 
We consider a sequence of birational maps 
$$
\overline{B_{0}^{l}}\dasharrow \overline{B_{1}^{l}}\dasharrow 
\cdots. 
$$ 
By \cite[Lemma 2.11]{matsuki}, 
after finitely many flips, 
$B_{i}^{l}$ does not contain $2$-dimensional 
flipping locus. 
Thus, we can assume that $B_{i}^{l}$ does not 
contain $2$-dimensional flipping locus for 
every $i\geq 0$. 
In particular, $B_i^{l}\dasharrow B_{i+1}^{l}$ is 
an isomorphism in codimension one for every $i\geq 0$. 

Next, we see $d_{S,b_{l-1}}(X,B)$. 
By the same argument as above, $B_{i}^{l-1}$ does not contain 
$2$-dimensional flipped locus 
after finitely many flips. 
Note that if $E$ is essential over $X_{i+1}$ and 
$a(E,X_{i+1}, B_{i+1})<1-b_{l-1}$, 
then $E$ is essential over $X_i$ and 
$a(E,X_i, B_i)\leq a(E,X_{i+1}, B_{i+1})$. 
It is because $B_i^l\dashrightarrow B_{i+1}^l$ 
is an isomorphism in codimension one. 
By using \cite[Lemma 2.11]{matsuki}, 
we see that $B_i^{l-1}$ does not contain $2$-dimensional 
flipping locus after finitely many steps.  
Therefore, we can assume that $B_i^{l-1}\dasharrow B_{i+1}^{l-1}$ is 
an isomorphism in codimension one for every $i\geq 0$

By repeating this argument, we can assume that 
$B_i^{j}\dashrightarrow B_{i+1}^{j}$ is an isomorphism 
in codimenison one for every $i$, $j$. 

If the log flip is of type $(1,2)$ or $(2,2)$, then 
$d_{S,b_0}(X,B)=d_{S, 0}(X,B)$ decreases by Lemma \ref{keisan}. 
Therefore, we can assume that all the flips are of type $(2,1)$ 
after finitely many steps. 
This sequence terminates by Lemma \ref{hom}.
\end{proof}

\begin{rem}[$3$-fold case] 
By using the weighted version of difficulty, 
we can easily prove the termination of $3$-fold 
canonical flips without using $\mathbb Q$-factoriality. 
This is \cite[4.14 Remark]{FA}. 
\end{rem}

\begin{ack} 
I was partially supported by the Inoue Foundation for 
Science. I talked about this paper at RIMS in May 2002. 
I would like to thank Professors Shigefumi Mori, Shigeru Mukai, 
and Noboru Nakayama, who gave me some comments. 
I am also grateful to Mr.~Kazushi Ueda, who answered my question about 
Borel-Moore homology. 
\end{ack}

\ifx\undefined\bysame
\newcommand{\bysame|{leavemode\hbox to3em{\hrulefill}\,}
\fi

\end{document}